\documentclass[11pt]{article}
\usepackage{newlfont,amsfonts,amssymb, oldgerm,amsmath,amsthm,amsgen,amscd,
}
\usepackage{epic}
\usepackage{eepic}
\usepackage{rotating}
\usepackage[dvips]{epsfig}
\begin{document}
\swapnumbers
\theoremstyle{definition}
\newtheorem{de}{Definition}[section]
\newtheorem{bem}[de]{Remark}
\newtheorem{bez}[de]{Notation}
\newtheorem{bsp}[de]{Example}
\theoremstyle{plain}
\newtheorem{lem}[de]{Lemma}
\newtheorem{satz}[de]{Proposition}
\newtheorem{folg}[de]{Corollary}
\newtheorem{theo}[de]{Theorem}

\newcommand{\bde}{\begin{de}}
\newcommand{\ede}{\end{de}}
\newcommand{\ul}{\underline}
\newcommand{\ol}{\overline}
\newcommand{\tbf}{\textbf}
\newcommand{\mc}{\mathcal}
\newcommand{\mb}{\mathbb}
\newcommand{\mf}{\mathfrak}
\newcommand{\bs}{\begin{satz}}
\newcommand{\es}{\end{satz}}
\newcommand{\btheo}{\begin{theo}}
\newcommand{\etheo}{\end{theo}}
\newcommand{\bfolg}{\begin{folg}}
\newcommand{\efolg}{\end{folg}}
\newcommand{\blem}{\begin{lem}}
\newcommand{\elem}{\end{lem}}
\newcommand{\bprf}{\begin{proof}}
\newcommand{\eprf}{\end{proof}}
\newcommand{\bd}{\begin{displaymath}}
\newcommand{\ed}{\end{displaymath}}
\newcommand{\be}{\begin{eqnarray*}}
\newcommand{\ee}{\end{eqnarray*}}
\newcommand{\eeqa}{\end{eqnarray}}
\newcommand{\beqa}{\begin{eqnarray}}
\newcommand{\bi}{\begin{itemize}}
\newcommand{\ei}{\end{itemize}}
\newcommand{\bnum}{\begin{enumerate}}
\newcommand{\enum}{\end{enumerate}}
\newcommand{\la}{\langle}
\newcommand{\ra}{\rangle}
\newcommand{\ve}{\varepsilon}
\newcommand{\vp}{\varphi}
\newcommand{\lra}{\longrightarrow}
\newcommand{\Lra}{\Leftrightarrow}
\newcommand{\Ra}{\Rightarrow}
\newcommand{\sub}{\subset}
\newcommand{\ems}{\emptyset}
\newcommand{\sms}{\setminus}
\newcommand{\ints}{\int\limits}
\newcommand{\sums}{\sum\limits}
\newcommand{\lims}{\lim\limits}
\newcommand{\bcup}{\bigcup\limits}
\newcommand{\bcap}{\bigcap\limits}
\newcommand{\beq}{\begin{equation}}
\newcommand{\eeq}{\end{equation}}
\newcommand{\einhalb}{\frac{1}{2}}
\newcommand{\rr}{\mathbb{R}}
\newcommand{\rn}{\mathbb{R}^n}
\newcommand{\ccc}{\mathbb{C}}
\newcommand{\cn}{\mathbb{C}^n}
\newcommand{\M}{{\cal M}}
\newcommand{\drehgleich}{\mbox{\begin{rotate}{90}$=$  \end{rotate}}}
\newcommand{\turngleich}{\mbox{\begin{turn}{90}$=$  \end{turn}}}
\newcommand{\turnsimeq}{\mbox{\begin{turn}{270}$\simeq$  \end{turn}}}
\newcommand{\vf}{\varphi}
\newcommand{\earr}{\end{array}\]}
\newcommand{\barr}{\[\begin{array}}
\newcommand{\bvec}{\left(\begin{array}{c}}
\newcommand{\evec}{\end{array}\right)}
\newcommand{\sumk}{\sum_{k=1}^n}
\newcommand{\sumi}{\sum_{i=1}^n}
\newcommand{\suml}{\sum_{l=1}^n}
\newcommand{\sumj}{\sum_{j=1}^n}
\newcommand{\suminf}{\sum_{k=0}^\infty}
\newcommand{\inv}{\frac{1}}
\newcommand{\wzbw}{\hfill $\Box$\\[0.2cm]}
\newcommand{\lag}{\mathfrak{g}}
\newcommand{\+}{\oplus}
\newcommand{\x}{\times}
\newcommand{\lx}{\ltimes}
\newcommand{\rrn}{\mathbb{R}^n}
\newcommand{\laso}{\mathfrak{so}}
\newcommand{\lason}{\mathfrak{so}(n)}
\newcommand{\w}{\omega}
\newcommand{\pmh}{{\cal P}(M,h)}
\newcommand{\s}{\sigma}
\newcommand{\deri}{\frac{\partial}}
\newcommand{\ddx}{\frac{\partial}{\partial x}}
\newcommand{\ddz}{\frac{\partial}{\partial z}}
\newcommand{\ddi}{\frac{\partial}{\partial y_i}}
\newcommand{\ddk}{\frac{\partial}{\partial y_k}}
\newcommand{\xz}{^{(x,z)}}
\newcommand{\mh}{(M,h)}
\newcommand{\wxz}{W_{(x,z)}}
\newcommand{\qmh}{{\cal Q}(M,h)}
\newcommand{\bbem}{\begin{bem}}
\newcommand{\ebem}{\end{bem}}
\newcommand{\bbez}{\begin{bez}}
\newcommand{\ebez}{\end{bez}}
\newcommand{\pr}{pr_{\lason}}
\newcommand{\huts}{\hat{\s}}
\newcommand{\whut}{\w^{\huts}}
\newcommand{\bhg}{{\cal B}_H(\lag)}
\newcommand{\aaa}{\alpha}
\newcommand{\bb}{\beta}
\newcommand{\lam}{\lambda}
\newcommand{\LL}{\Lambda}
\newcommand{\D}{\Delta}
\newcommand{\ß}{\beta}
\newcommand{\ä}{\alpha}
\newcommand{\W}{\Omega}
\newcommand{\esel}{\ensuremath{\mathfrak{sl}(2,\ccc)}}

\bibliographystyle{alpha}


\title{Towards a classification of Lorentzian holonomy groups. Part II: Semisimple, non-simple weak-Berger algebras}

\author{Thomas Leistner
}

\maketitle

\begin{abstract}
The holonomy group of an $(n+2)$--dimensional simply-connected, indecomposable but non-irreducible
Lorentzian manifold $(M,h)$ is contained in
the parabolic group $( \mathbb{R} \times
SO(n) )\ltimes \mathbb{R}^n$. The main ingredient of such a holonomy group is
the $SO(n)$--projection $G:=pr_{SO(n)}(Hol_p(M,h))$ and one may ask whether it has to be
a Riemannian holonomy group. In this paper we show that this is always the case, completing our
results of \cite{leistner03}. We draw consequences for the existence of parallel spinors on Lorentzian manifolds.
\end{abstract}
\setcounter{tocdepth}{1}
\tableofcontents

\section{Introduction}

This paper is an addendum to our paper \cite{leistner03} where we gave a partial
classification of reduced holonomy groups of indecomposable Lorentzian
manifolds. The holonomy group of an $(n+2)$-dimensional, indecomposable, non-irreducible Lorentzian
manifold is contained in the parabolic
group whose Lie algebra is $(\rr \+ \lason)\ltimes \rrn$. Concerning the three projections, L. Berard-Bergery and
A. Ikemakhen distinguished  in \cite{bb-ike93}
four different types of indecomposable subalgebras
of $(\rr \+ \lason)\ltimes \rrn$.
But the main ingredient of such a holonomy algebra is the $\lason$-projection.
Also in \cite{bb-ike93} a Borel-Lichnerowicz-type decomposition property is
proved:
\btheo\label{theoII}\cite{bb-ike93}
Let be $\lag=pr_{\lason}\mf{hol}(M^{n+2},h)$.
Then it holds: For the decomposition of the representation space $\rrn=E_0 \+\ldots
\+ E_k$ of $\lag$ into irreducible components there is
a decomposition of the Lie algebra $\lag=\lag_1\+\ldots\+\lag_k \subset \lason$ into ideals
and each of these
$\lag_i$ acts  irreducibly on $E_i$ and trivial on $E_j$ for $i\not=j$.
\etheo
Thus for
 the problem of classifying possible $\lason$-parts of
indecomposable Lorentzian holonomy algebras, one can restrict oneself to the study of irreducibly acting
subalgebras of $\lason$. For this we introduced in \cite{leistner02} the notion of a weak-Berger algebra and showed that
the $\lason$-component of an indecomposable Lorentzian holonomy algebra is a
weak-Berger algebra and all its irreducibly acting components too. Using this we classified the unitary
acting weak-Berger algebras.
In the first part \cite{leistner03} to the present paper we extended the classification to
irreducible weak-Berger algebras which are simple, obtaining the following
theorem:

\btheo
\cite[Theorem 3.2]{leistner02}, \cite[Theorem 3.21 and Theorem  4.7]{leistner03}
Let $\lag\subset  \laso(n,\rr)$ be a real, irreducible weak-Berger algebra which is simple or acts unitary.
Then it is a Riemannian holonomy algebra.
\etheo

For this we used the distinction of real representations into representations of real
and of non-real type. Orthogonal representation of non-real type are unitary
representations. For representations of real type which are weak-Berger
the complexified Lie algebra (with the complexified representation of
course) is also weak-Berger. Furthermore holds that an orthogonal Lie algebra
of real type has to be semisimple. Hence \cite{leistner03} leaves open the
following problem: Classify all complex, irreducible weak-Berger
algebras which are semisimple but not simple. We will solve this problem  in the
present paper. It uses very much results and
proofs of \cite{leistner03}. Here we will complete the proof of the
following theorem by completing the classification result in the semisimple
case.
\btheo\label{theo1}
Any weak-Berger algebra is the holonomy algebra of a Riemannian holonomy
algebra. In particular: The $\lason$-component of an indecomposable, non-
irreducible Lorentzian holonomy algebra is a Riemannian holonomy algebra.
\etheo
To complete the proof in the semisimple case we will show the following:
\bs \label{satz1}
Any irreducibly acting, semisimple, non-simple complex weak-Berger algebra is
the compexification of Riemannian holonomy algebra.
\es
To prove this statement we will use several results of \cite{leistner03} which
describe the weak-Berger property in terms of root systems.

Up to dimension $n=9$ the result of theorem \ref{theo1} was proved by A. Galaev in \cite{galaev}, partially using
results of \cite{leistner02}.

At the end of this paper we will show some applications for Lorentzian manifolds with parallel spinors showing that
their holonomy group is of the form $G\ltimes \rrn$ where $G$ is a holonomy group of a Riemannian manifold with
parallel spinor.

\section{Proof of the result}
\subsection{The notion of a weak-Berger algebra}
First we recall the notion of a weak-Berger algebra.
\begin{de}\label{weakdef}
Let ${\frak g}\subset \mathfrak{so}( E,h)$ be an orthogonal Lie algebra. Then ${\frak g}$ is
is called {\bf weak-Berger algebra} if $\lag=span\{ Q(x) |  x\in E, Q\in {\cal B}_h({\frak g})\}$, where
${\cal B}_h({\frak g}):= \{ Q\in E^* \otimes {\frak g}  | h(Q(x)y,z) + h(Q(y)z,x)
+ h(Q(z)x,y)=0\}$.
\end{de}
If $E$ is a real vector space and the the real Lie algebra $\lag_0\subset \laso(E,h)$ is irreducible of real type, i.e.
the complexification is irreducible,  then
$\lag_0$ is a weak-Berger algebra if and only if $\lag:=\lag_0^\ccc\subset \laso(E^\ccc, h^\ccc)$ is a
complex weak-Berger algebra.
Since an irreducibly acting, complex $\lag\subset \laso(V,H))$ --- for $V$ a complex vector space ---
is semisimple we can use the tools of root space and weight space decomposition to
classify complex weak-Berger algebras.

We denote by $\Omega$ the weights of $\lag\subset \laso(V,H)$,
by $\Delta$ the roots of $\lag$ and we set $\Delta_0=\Delta\cup\{0\}$.
If $\alpha\in \Delta$ the we define $\Omega_\alpha:=
\{\mu\in\Omega|\mu+\alpha\in \Omega\}$. Then it holds:

\bs\label{musatz}\cite[Proposition 2.6]{leistner03}
Let $\lag$ be a semisimple
Lie algebra with roots $\Delta$ and $\Delta_0=\Delta\cup\{0\}$. Let $\lag\subset
\laso(V,H)$ irreducible, weak-Berger with weights $\Omega$. Then the following properties are satisfied:
\begin{description}
\item[(PI)] There is a $\mu\in \Omega$ and a hyperplane $U\subset \mf{t}^*$ such that
\beq \label{mu}
\Omega\subset \left\{ \mu+\beta\ |\ \beta\in \Delta_0\right\}\cup U\cup
\left\{- \mu+\beta\ |\ \beta\in \Delta_0\right\}.
\eeq
\item[(PII)] For every $\alpha\in \Delta$ there is a $\mu_\alpha\in \Omega$ such that
\beq \label{mua}
\Omega_\alpha\subset \left\{ \mu_\alpha-\alpha+\beta\ |\ \beta\in \Delta_0\right\}\cup
\left\{ -\mu_\alpha+\beta\ |\ \beta\in \Delta_0\right\}.
\eeq
\end{description}
\es
Furthermore we obtained:
\bs\label{triplesatz}\cite[Proposition 2.13]{leistner03}
Let $\lag\subset \laso(V,H)$ be an irreducible complex weak-Berger algebra. Then there is an extremal weight
$\Lambda$ such that
one of the following properties is satisfied:
\begin{description}
\item[(SI)]
There is a pair $(\Lambda, U)$ with an extremal weight $\Lambda$   a hyperplane $U$
in $\mf{t}^*$ such that every extremal weight different from $\LL$ and $-\LL$ is
contained in $U$ and
$\Omega\subset
 \left\{ \LL+\beta\ |\ \beta\in \Delta_0\right\}\cup U
 \cup
\left\{ \LL+\beta\ |\ \beta\in \Delta_0\right\}$.
\item[(SII)] There is an $\ä\in \Delta$ such that
$
\Omega_\alpha\subset \left\{ \Lambda-\alpha+\beta\ |\ \beta\in \Delta_0\right\}\cup
\left\{ -\Lambda+\beta\ |\ \beta\in \Delta_0\right\}.
$
\end{description}
There is a fundamental system such that the extremal weight in (SI) and (SII) is the highest weight.
\es
In the following we classified all simple Lie algebras which satisfy (PI) and (PII) or (SI) or (SII) and showed that they are
Riemannian holonomy algebras. In the present paper we will do this for semisimple Lie algebras which are not simple.

\subsection{Semisimple, non-simple weak-Berger algebras}
From now on let $\lag$ be a complex semisimple, non-simple Lie algebra, irreducible represented on a complex vector
space $V$. To a decomposition of $\lag$ into ideals $\lag=\lag_1\+\lag_2$ corresponds a decomposition of
the irreducible module $V$ into factors $V=V_1\otimes V_2$ which are irreducible $\lag_1$- resp. $\lag_2$-modules.
$X=(X_1,X_2)\in \lag$ acts as follows: $X\cdot (v_1\otimes v_2)= (X_1\cdot v_1)\otimes v_2 + v_1\otimes X_2\otimes v_2$.
The Cartan subalgebra $\mf{t}$ of $\lag$ is the sum of the Cartan subalgebras of $\lag_1$ and $\lag_2$. If
$\Delta $ are the roots of $\lag$ and $\Delta^i$ the roots of $\lag_i$ then $\Delta=\Delta^1\cup\Delta^2$.
For the weights it holds $\Omega=\Omega^1+\Omega^2$. Analogously we denote for $\aaa\in \Delta^i$ the set
$\Omega_\aaa^i$.
Then holds the following
\blem\label{lemma1}
Let $\lag=\lag_1 \+\lag_2$ be a semisimple Lie algebra, $V=V_1 \otimes V_2$ be an irreducible representation of it.
If $\aaa\in \Delta^1$ then it holds
\beq\label{omegaalpha}
\W_\aaa=\W^1_\aaa +\W^2
\eeq
\elem
\bprf
For $\lam\in \W_\aaa$ we have $\W\ni \lam+\aaa=\lam_1+\aaa+\lam_2$ with $\lam_i\in \W_i$. Hence
$\lam_1+\aaa\in \W^1$. If otherwise $\lam_1+\aaa\in \W^1$ then $\lam_1+\lam_2+\aaa\in \W$, i.e.
$\lam_1+\lam_2\in \W_\aaa$.
\eprf
Assuming the weak-Berger property we get by this:
\blem\label{lemma2}
Let $\lag=\lag_1 \+\lag_2$ be a semisimple Lie algebra, $V=V_1 \otimes V_2$ be an irreducible representation of it which
is weak-Berger. If the dimensions of $V_1$ and $V_2$ are greater than $2$, then for any $\aaa\in \Delta^i$ the set
$\W^i_\aaa$ contains at most 2 elements.
\elem
\bprf
Suppose that $dim\ V_2\ge 3$, i.e. $\#\W^2\ge 3$. Let $\aaa\in \Delta^1,\ \lam_1\in \W^1_\aaa$ and
$\lam_2\in\W^2$, i.e. $\lam_1+\lam_2\in \W^1_\aaa$.
Now from the property (PII) follows that there is a $\mu_\aaa=:\mu^1_\aaa+\mu^2_\aaa \in \W$
such that
$\lam_1+\lam_2=\mu_\aaa-\aaa+\beta$ or $\lam_1+\lam_2=-\mu_\aaa+\beta$ with $\beta\in \Delta_0=\Delta^1\cup\Delta^2\cup
\{0\}$.
If now $\# \Omega^2\ge 3$ and $\#\W^1_\aaa\ge 3$ then we can choose $\lam_1\not= \mu^1_\aaa-\aaa$, $\lam_1\not=
-\mu^1_\aaa$ and $\lam_2\not=\pm\mu^2_\aaa$. This gives a contradiction.
\eprf
Now we can use a result of L. Schwachh\"{o}fer from \cite{schwachhoefer2}.
\bs \cite[Lemma 3.23]{schwachhoefer2}\label{satzschwachhoefer}
Let $\lag\subset\mf{gl}(n,\ccc)$ be an irreducibly acting, semisimple subalgebra. If for any $\aaa$ the set
$\W_\aaa$ contains at most two elements, then $\lag$ is conjugate to one of the following representations:
\bnum
\item
$\mf{sl}(n,\ccc)$ acting on $\ccc^n$; in this case $\W_\aaa$ is a singleton for all $\aaa\in \Delta$.
\item
$\mf{so}(n,\ccc)$ acting on $\ccc^n$; in this case $ \W_\aaa$ contains two elements for all $\aaa\in \Delta$,
and their sum equals to $-\aaa$.
\item
$\mf{sp}(n,\ccc)$ acting on $\ccc^{2n}$; in this case $ \W_\aaa$ contains two elements  $\aaa\in \Delta$ is short,
and their sum equals to $-\aaa$, and $\W_\aaa=\{-\einhalb \aaa\}$.
\enum
\es
From this result we obtain the following corollary, proving proposition \ref{satz1} if the dimensions of the factors of $V$
are greater than 2.
\bfolg\label{folg1}
Let $\lag\subset \mf{so}(V,h)$ be a complex, semisimple, non-simple, irreducibly acting weak-Berger algebra. If $\lag$ decomposes into
$\lag=\lag_1\+\lag_2$ such that for the corresponding decomposition of $V=V_1\otimes V_2$ holds that
$dim\ V_i\ge 3$ for $i=1,2$, then it holds:
$\lag=\mf{so}(n,\ccc)\+\mf{so}(m,\ccc)$ acting on $\ccc^n\otimes \ccc^m$, or
$\lag=\mf{sp}(n,\ccc)\+\mf{sp}(m,\ccc)$ acting on $\ccc^{2n}\otimes \ccc^{2m}$. In particular it is the complexification of a
Riemannian holonomy representation of a symmetric space of type BDI resp. CII (for the types see \cite{helgason78}.
\efolg
\bprf
By lemma \ref{lemma2} it must hold $\#\W^i_\aaa\le 2$ for both summands. So
we have to built sums of the Lie algebras of
proposition \ref{satzschwachhoefer}. But only the sum  of two orthogonal acting Lie algebras, or a sum of
two symplectic acting Lie algebras acts orthogonal.
\eprf
By this result we are left with semisimple Lie algebras where the irreducible representation of one summand is two-dimensional, i.e. $\lag=\mf{sl}(2,\ccc)\+\lag_2$ and $V=\ccc^2\otimes V_2$.
Since we are interested in $\lag\subset \mf{so}(V,h)$ and $\mf{sl}(2,\ccc)$
acts symplectic on $\ccc^2$ the representation of $\lag_2$ on $V_2$ has to be symplectic too.

In this situation we prove the following
\bs \label{satz2}
Let $\lag=\mf{sl}(2,\ccc)\+\lag_2$ be a semismple, complex weak-Berger algebra, acting irreducibly on
$\ccc^2\otimes V_2$. Then $\lag_2\subset \mf{sp}(V_2)$ satisfies the following properties:
\begin{description}
\item[(PIII)] There is a $\mu\in \Omega^2$ and an affine hyperplane $A\subset \mf{t}_2^*$ such that
\beq \label{muIII}
\Omega^2\subset \left\{ \mu+\beta\ |\ \beta\in \Delta^2_0\right\}\cup A\cup
\left\{- \mu\right\}.
\eeq
\item[(PIV)]
There is a $\mu\in \Omega^2$ such that
\beq \label{muIV}
\Omega^2\subset \left\{ \mu+\beta\ |\ \beta\in \Delta^2_0\right\}\cup
\left\{- \mu+\beta\ |\ \beta\in \Delta^2_0\right\}.
\eeq
\end{description}
\es
\bprf
Since $\lag$ is weak-Berger it satisfies the properties (PI) and (PII). We draw the consequences from both for
$\lag=\mf{sl}(2,\ccc)\+\lag_2$. For the representation of $\mf{sl}(2,\ccc)$ on $\ccc^2$ we have that
$\W^1=\{\Lambda,-\Lambda\}$. Let $\alpha\in \Delta^1=\{\alpha,-\alpha\}$ be
the positive root of $\mf{sl}(2,\ccc)$. Hence $\W^1_ \aaa=\{-\Lambda\}$, since $-\Lambda=\Lambda-\aaa$.
\begin{description}
\item{(PIII)} $\lag$ satisfies the property (PI) with a hyperplane $U:=(T_1+T_2)^\bot$ and a weight $\mu=\mu_1+\mu_2$.
Let $\lam=\mu_1+\lam_2\in \W$ be a weight of $\lag$.
If $\lam$ lies in a hyperplane of
$\mf{t}=\mf{t}_1\+^\bot \mf{t}_2$, then $0=\la \mu_1,T_1\ra+\la\lam_2,T_2\ra$, i.e.
$\lam_2$ lies in an affine hyperplane of $\mf{t}$.
If $\lam=\mu+\beta=\mu_1+\mu_2+\beta$,  then $\lam_2=\mu_2+\beta$ with $\beta\in \Delta^2$.
If $\lam=-\mu+\beta=-\mu_1-\mu_2+\beta$, then $\beta$ has to be in $\Delta^1$ and $\lam_2=\mu_2$. Hence
$\lag_2\subset \mf{sp}(V_2) $ satisfies (PIII).
\item{(PIV)} $\lag$ satisfies the property (PII). Suppose that $\aaa$ is the positive root of $\mf{sl}(2,\ccc)$.
Then $\W^1_\aaa=\{-\Lambda\}$ and $\W_\aaa=\{-\Lambda\}\cup\W^2$. Now let $\lam\in \W^2$, i.e.
$-\Lambda+\lam\in \W_\aaa$. By (PII) there is a $\mu_\aaa=\mu^1+\mu^2$ such that
$-\Lambda+\lam=\mu_\aaa-\aaa+\beta$ or $-\Lambda+\lam=-\mu_\aaa+\beta$ with $\beta\in
\Delta_0=\Delta^1\cup\Delta^2\cup\{0\}$. Since $\mu_1=\pm\Lambda$ this implies
$\lam\in \{\mu_2 +\beta|\beta\in \Delta_0^2\}\cup \{-\mu_2 +\beta|\beta\in \Delta_0^2\}$, i.e. (PIV) is satisfied.
\end{description}
So we have shown that both, (PIII) and (PIV) are satisfied.
\eprf

\begin{bsp}
We set $\lag_2=\esel$ and check  if $\lag=\esel\+\esel$ acting on $\ccc^2\otimes V_2$ is a weak-Berger algebra.
This is to check whether $\esel$ acting on $V_2$ satisfies (PIII) and (PIV) and is symplectic.
To be symplectic means that the representations has an
even number of weights, (PIV) implies that $V_2$ has at most 6 weights but (PIII) implies that $V_2$ has at most 4 weights.
Hence the only weak-Berger algebras with the structure of $\esel\+\esel$ are those acting on $\ccc^4$ and on $\ccc^2\otimes
sym_0^3\ccc^2=\ccc^8$. Both are of course complexifications of Riemannian holonomy representations, the first of the
4-dimensional symmetric space of type CII, i.e. $Sp(2)/Sp(1)\cdot SP(1)$ and the second of the
8-dimensional symmetric space of type GI, i.e. $G_2/SU(2)\cdot SU(2)$ (in the compact case, see \cite{helgason78}).
\end{bsp}

Now we try to reduce the problem in a way that we only have to deal with simple Lie algebras.
\blem\label{lemma3}
Let $\lag\subset\mf{gl}(V)$ be a semisimple, complex Lie algebra acting irreducibly on $V$, satisfying the property
(PIV).
Then $\lag$ is simple or $\lag=\esel\+\lag_2$ acting on $\ccc^2\otimes V_2$.
\elem
\bprf
Suppose that $\lag=\lag_1\+\lag_2$ and that $\#\Omega^1\ge 3$. Let $\mu=\mu_1+\mu_2$ be the weight from the
property (PIV). We consider a weight $\lam=\lam_1+\lam_2\in\W=\W^1+\W^2$ with $\lam_1\not=\pm\mu_1$. Then (PIV) implies that
$\lam_2=\mu_2$ or $\lam_2=-\mu_2$, i.e. $\#\W^2\le 2$. This implies the proposition of the lemma.
\eprf

To complete the reduction we need a further
\blem\label{lemma4}
Let $\lag=\esel\+\esel\+\lag_3$ be a semisimple complex Lie algebra, acting irreducibly on
$\ccc^2\otimes\ccc^2\otimes V_3$ and
satisfying the property (PII). Then for any root $\aaa\in \Delta^3$ of $\lag_3$ holds
$\#\W^3_\aaa\le 2$.
\elem
\bprf
Let $\aaa\in \Delta^3$ and $\mu_\aaa^1+\mu_\aaa^2+\mu_\aaa^3$ the weight from the property (PII). Then
$\W_\aaa=\W^1+\W^2+\W^3_\aaa\ni -\mu_\aaa^1+\mu_\aaa^2+\lam$ with $\lam\in \W^3_\aaa$ arbitrary.
Again (PII) implies $\lam=\mu_\aaa^3-\aaa$ or $\lam=-\mu_\aaa^3$, i.e. $\#\W^3_\aaa\le 2$.
\eprf
Both lemmata give the following result.
\bs
Let $\lag=\esel+ \lag_2$ be a semisimple, complex Lie algebra, acting irreducibly on $\ccc^2\otimes V_2$
which is supposed to be weak-Berger.
Then $\lag_2$ is simple, acts irreducible and symplectic on $V_2$ satisfying (PIII) and (PIV), or
$\lag=\esel\+\esel\+\laso(n,\ccc)=\laso(4,\ccc)\+\laso(n,\ccc)$ acting irreducibly on $\ccc^4\otimes \ccc^n$.
\es
\bprf
The proof is obvious by lemma \ref{lemma3} and lemma \ref{lemma4} and the result of proposition \ref{satzschwachhoefer}
keeping in mind that $\lag$ is orthogonal, hence $\lag_2$ is symplectic and $\lag_3$ has to be orthogonal again.
\eprf
Of course, the representation of $\laso(4,\ccc)\+\laso(n,\ccc)$  on $\ccc^4\otimes \ccc^n$
is the complexification of a Riemannian holonomy representation of the symmetric space of type BDI.

\subsection{Simple Lie algebras satisfying (PIII) and (PIV)}

In this section we deal with the remaining problem to classify complex, simple irreducibly acting symplectic Lie algebras
with the property (PIII) and (PIV).

\bs
Let $\lag\subset \mf{sp}(V)$ be simple, irreducibly acting and satisfying (PIV).
Then it satisfies  (SII).
\es
\bprf
First we note that the fact that the representation is symplectic leaves us with the simple Lie algebras with root systems $A_n, B_n, C_n, D_n$ and $E_7$.
In particular the Lie algebra of type $G_2$ is excluded. This
implies that for two roots $\aaa$ and $\beta$ it holds that
$\left|\frac{\la\aaa,\beta\ra}{\|\aaa\|^2}\right|\in \{1,\einhalb,0\}$, a fact which we will use several times in the following proof.

Let $\mu$ be the weight from the property (PIV). We consider two cases.
\begin{description}
\item{{\em Case 1: $\mu$ is not an extremal weight:}}
In this case there is a root $\aaa\in \Delta$ such that $\mu+\aaa=\Lambda$ is extremal. We show that (SII) is satisfied
with the tripel $(\Lambda,-\Lambda,\aaa)$.

We suppose that (SII) is not satisfied, i.e. there is a $\lam\in \W_\aaa\subset\W$
such that neither $\lam=\LL-\aaa+\beta$ nor $\lam=-\LL+\beta$ for a $\beta\in \Delta$.
$\lam\in\W$ and $\lam+\aaa\in \W$ gives by (PIV) that
 $\lam=-\LL+\aaa+\beta$ with $\beta\in \Delta$ and $\aaa+\beta\not\in\Delta_0$, as well as
$\lam=\LL-2\aaa+\gamma$ with $\gamma\in \Delta$ and $\aaa-\gamma\not\in\Delta_0$.
By properties of root systems this implies that $\la\aaa,\beta\ra \ge 0$ and $\la\aaa,\gamma\ra\le 0$.
Furthermore it is
\beq
\label{2ll1}
2\LL=3\aaa+\beta-\gamma.
\eeq
Now it is
$\frac{2\la\LL,\aaa\ra}{\|\aaa\|^2}= 3+ \frac{\la\beta,\aaa\ra}{\|\aaa\|^2}-\frac{\la\gamma,\aaa\ra}{\|\aaa\|^2}\ge 3$,
entailing
$\LL-3\aaa\in\W$. Since $\LL-3\aaa\not=\LL-\aaa+\delta$ for a $\delta\in \Delta_0$. (PIV) implies
$\LL-3\aaa\not=-\LL+\aaa+\delta$, i.e.
\beq
\label{2ll2}
2\LL=4\aaa+\delta,
\eeq
with $\delta\not=-\aaa$.
(\ref{2ll1}) and (\ref{2ll2}) give
\beq
\label{2ll3}
0=\aaa+\delta+\gamma-\beta.
\eeq

Now suppose that $\frac{2\la\LL,\aaa\ra}{\|\aaa\|^2}= 3$, i.e. $\la\beta,\aaa\ra=\la\gamma,\aaa\ra=0$.
In this case (\ref{2ll2}) gives
$\frac{2\la\delta,\aaa\ra}{\|\delta\|^2}=-2$ and therefore $\frac{2\la\LL,\delta\ra}{\|\delta\|^2}= -3$.
This implies that $\LL+3\delta\in \W$, but this is together with (\ref{2ll2}) is a contradiction to (PIV).

Now suppose that $\frac{2\la\LL,\aaa\ra}{\|\aaa\|^2}= 4$, i.e.
$\frac{\la\beta,\aaa\ra-\la\gamma,\aaa\ra}{\|\aaa\|^2}=1$. Then (\ref{2ll3})  implies
$\la \aaa,\delta\ra=0$. $\LL-4\aaa\in \W$ implies by (PIV) and $3\aaa\not\in\Delta$ that
$2\LL=5\aaa+\ve$, i.e. $\aaa-\delta\in \delta$. Since $\la\aaa,\delta\ra=0$ this implies that
$\aaa$ and $\delta$ are short roots and $\aaa-\delta $ is a long one, i.e.
$\frac{\|\delta\|^2}{\|\aaa-\delta\|^2}=\frac{\|\aaa\|^2}{\|\aaa-\delta\|^2}=\einhalb$. But this gives
that
$\frac{2\la\LL,\aaa-\delta\ra}{\|\aaa-\delta\|^2}=\frac{5}{2}$ which is a contradiction.

Finally suppose that
$\frac{2\la\LL,\aaa\ra}{\|\aaa\|^2}\ge 5$.
Hence $\frac{\la\beta,\aaa\ra}{\|\aaa\|^2}-\frac{\la\gamma,\aaa\ra}{\|\aaa\|^2}\ge 2$
On the other hand $\LL-5\aaa\in \W$ and by (PIV)
$2\aaa-\delta\in \Delta$. This implies that
$\frac{2\la\aaa,\delta\ra}{\|\aaa\|^2}\ge 2$. But both inequalities are a contradiction to
(\ref{2ll3}).
\item{{\em Case 2. $\mu:=\LL$ is an extremal weight.}} To proceed analogously as in the first case
we fix a root $\aaa\in\D$, which is supposed to be long in case of root systems with roots of different lenght,
 and we show that (SII) is satisfied for the tripel $(\LL,-\LL,\aaa)$.

Again we suppose that (SII) is not satisfied, i.e. there is a $\lam\in \W_\aaa\subset\W$
such that neither $\lam=\LL-\aaa+\beta$ nor $\lam=-\LL+\beta$ for a $\beta\in \Delta$.
$\lam\in\W$ and $\lam+\aaa\in \W$ gives by (PIV) that
 $\lam=\LL+\beta$ with $\beta\in \Delta$ and $\aaa+\beta\not\in\Delta_0$, as well as
$\lam=-\LL-\aaa+\gamma$ with $\gamma\in \Delta$ and $\aaa-\gamma\not\in\Delta_0$.
By properties of root systems this implies that $\la\aaa,\beta\ra \ge 0$ and $\la\aaa,\gamma\ra\le 0$.
Since $\aaa$ is supposed to be a long root this the same as $\frac{\la \aaa,\beta\ra}{\|\aaa\|^2}\in\{0,\einhalb\}$ and
$\frac{\la \aaa,\gamma\ra}{\|\aaa\|^2}\in\{-\einhalb,0\}$.
Furthermore it is
\beq
\label{2ll4}
2\LL=-\aaa-\beta+\gamma
\eeq
and hence $\mathbb{Z}\ni\frac{2\la\LL,\aaa\ra}{\|\aaa\|^2}= 1- \frac{\la\beta,\aaa\ra}{\|\aaa\|^2}+\frac{\la\gamma,\aaa\ra}{\|\aaa\|^2}
=:a\le -1$. Then of course $a\in \{-2,-1\}$.

First suppose that $a=-1$. In the case it is $\la \aaa,\beta\ra=\la\aaa,\gamma\ra=0$. Then because of
$\mathbb{Z}\ni\frac{2\la\LL,\beta\ra}{\|\beta\|^2}=-1+ \frac{\la\beta,\gamma\ra}{\|\beta\|^2}$ and
$\mathbb{Z}\ni\frac{2\la\LL,\gamma\ra}{\|\gamma\|^2}=-1+ \frac{\la\beta,\gamma\ra}{\|\gamma\|^2}$ it must hold that
$\frac{\la\beta,\gamma\ra}{\|\beta\|^2}$ and $\frac{\la\beta,\gamma\ra}{\|\gamma\|^2}$ are integers. But this can only be true
if $\beta$ and $\gamma$ are both, long and short. This is impossible.

Now suppose that $a=- 2$, i.e.
$\frac{\la \aaa,\beta\ra}{\|\aaa\|^2}=\einhalb $ and
$\frac{\la \aaa,\gamma\ra}{\|\aaa\|^2}=-\einhalb$.
Then $\LL-2\aaa\in \W$, i.e. by (PIV) we get that
\beq
\label{2ll5}
2\LL=-2\aaa+\delta.
\eeq
with $\delta\in \Delta_0$ with $\delta\not=\pm \aaa$ because otherwise we would get $a=-1$ or $a=-3$.

Now the existence of a root $\ve$ with the property $\|\delta\|\le\|\ve\|$ would give a contradiction since
\[\mathbb{Z}\ni \frac{2\la\LL,\ve\ra}{\|\ve\|^2}=-\underbrace{\frac{2\la\aaa,\ve\ra}{\|\ve\|^2}}_{\in\mathbb{Z}}
+\underbrace{\frac{\la\delta,\ve\ra}{\|\ve\|^2}}_{\not\in\mathbb{Z}}.\]
This implies that $\delta$ is a long root in the root system of type $C_n$. In $C_n$ the system of long roots
equals to $A_1\times \ldots \times A_1$. By this $\frac{\la \aaa,\beta\ra}{\|\aaa\|^2}=\einhalb $ and
$\frac{\la \aaa,\gamma\ra}{\|\aaa\|^2}=-\einhalb$ implies that $\beta$ and $\gamma$ are short roots recalling that
$\aaa$ was supposed to be a long one.

But then by (\ref{2ll4}) we get
\[
\frac{2\la\LL,\beta\ra}{\|\beta\|^2}=
- \frac{\la\aaa,\beta\ra}{\|\beta\|^2}+\frac{\la\beta,\gamma\ra}{\|\beta\|^2}-1
=- \einhalb\frac{\|\aaa\|^2}{\|\beta\|^2}+\frac{\la\beta,\gamma\ra}{\|\beta\|^2}-1
=-2+\frac{\la\beta,\gamma\ra}{\|\beta\|^2}\not\in\mathbb{Z}\]
since $\beta$ and $\gamma$ are short. But this is a contradiction.

\end{description}
\eprf

As a consequence of this proposition we only have to check the irreducible representations of simple
Lie algebra whether they satisfy (SII) --- done in \cite{leistner03} --- and then to add the
condition that the representations are symplectic --- instead of orthogonal. We obtain the following result.

\bs
Let $\lag\subset\mf {sp}(V)$ be a complex, simple, irreducibly and symplectic acting Lie algebra satisfying (PIII) and
(PIV) and different from $\esel$.
Then the root system and the highest weight of the representation are one of the following:
\bnum
\item $A_5$: $\w_3$, i.e. $\lag=\mf{sl}(6,\ccc)$ acting on $\wedge^3\ccc^6$.
\item $C_n$: $\w_1$, i.e. $\lag=\mf{sp}(n,\ccc)$ acting on $\ccc^{2n}$.
\item $C_3$: $\w_3$, i.e. $\lag=\mf{sp}(3,\ccc)$ acting on $\ccc^{14}$.
\item $D_6$: $\w_6$, i.e. $\lag=\mf{so}(12,\ccc)$ acting on $\ccc^{32}$ as spinor representation.
\item $E_7$: $\w_1$, i.e. the standard representation of $E_7$ of dimension $56$.
\enum
\es
\bprf
(PIV) implies (SII), so we use former results checking
the Lie algebras satisfying (SII) whether they are symplectic. For this we consider two cases.

First we suppose that $0\in \W$. In proposition 3.6 and corollary 3.9 of \cite{leistner03}
it is proved that any such representation which satisfies (SII) and is self-dual is orthogonal.
Hence if $0$ is a weight, no symplectic representation satisfies (SII).

Now suppose that $0\not\in\W$. In the proof of proposition 3.18 of \cite{leistner03} we have shown that the
representations of the following Lie algebras with $0\not\in \W$ satisfy (SII). Now we check if these are symplectic
and in some cases if they satisfy (PIII) and (PIV).
\bnum
\item $A_n$ with $n\le 7$ odd, $\LL=\w_{\frac{n+1}{2}}$. The only representation of these which is symplectic is the
one for $n=5$.
\item $B_n$: $\w_n$ for $n\le 7$ the spin representations, and $\w_1+\w_2$ for $n=2$. The latter is symplectic and the
the former is symplectic for $n=5,6$. ($B_2\simeq C_2$ we will study in the next point.) Now we show that these remaining
representations does not satisfy (PIII) or (PIV).

Of course the representation of $B_2$ with highest weight $\LL=\w_1+\w_2=\frac{3}{2}e_1+\einhalb e_2$ can not satisfy
(PIV) because it has 12 weights while $B_2$ has only 8 roots.

The spin representation for $n=6$ can not obey (PIV): W.l.o.g we may assume that $\LL $ from (PIV) is the highest weight
$\LL=\einhalb(e_1+ \ldots +e_6)$. But then for the weight $\lam=\einhalb (e_1+e_2+e_3-e_4-e_5-e_6)$ it holds neither
$\LL-\lam\in \Delta_0$ nor $\LL+\lam\in \Delta_0$.

The spin representation for $n=5$ does satisfy (PIV) but  not (PIII)
since all the weights $\einhalb(\pm e_1\pm \ldots \pm e_5)$
with 3 minus signs can not lie on the same affine hyper plane.

Hence none of the symplectic representations satisfying (SII) satisfies (PIII) and (PIV).
\item $C_n$ with $\LL=\w_1+\w_i$ or $\LL=\w_i$. These are symplectic for $i$ even in the first case and
for $i$ odd in the second case.

Again we have to impose the condition (PIV) on both. First we consider the representation with
highest weight $\LL=\w_i=e_1+\ldots +e_i$. Since the set of roots of $C_n$ equals
to $\{e_i\pm e_j, \pm2e_k\}$ we get
\[\W=\{\pm e_{k_1} \pm\ldots \pm e_{k_i}\}\cup
\{\pm e_{k_1} \pm\ldots \pm e_{k_{i-2}}\}\cup
\ldots\cup \{\pm e_k\}.\]
From this one sees that (PIV) can not be satisfied if $n\ge 5$.

With analogous considerations we exclude the case where $\LL=\w_1+\w_i$ with $i$ even.
\item $D_n$ with $\LL=\w_n$ and $n\le 8$. But these are only symplectic for $n=6$ and $n=2$. The latter
is excluded since $D_2=A_1\x A_1$, a case which is handled in the previous subsection.
\item For $E_7$ remains only the representation given in the proposition.
\enum\eprf

If we now combine the result of this and the previous subsection we get the following:
\bfolg
Let $\lag=\esel\+\lag_2$ be a semisimple, complex weak-Berger algebra acting on $\ccc^2\otimes V_2$.
Then it is the complexification of a Riemannian holonomy representation, in particular the complexification of the
holonomy representation of a non-symmetric $Sp(1)\cdot Sp(n)$-manifold or of the
following Riemannian symmetric spaces (we list only the compact symmetric space):
\bnum
\item Type $EII$: $E_6/SU(2)\cdot SU(6)$,
\item Type $CII$: $Sp(n+1)/Sp(1)\cdot Sp(n)$,
\item Type $FI$: $F_4/SU(2)\cdot Sp(3)$,
\item Type $EVI$: $E_7/ SU(2)\cdot Spin(12)$,
\item Type $EIX$: $E_8/SU(2)\cdot E_7$
\enum
and of type $GI$, i.e. $G_2/SU(2)\cdot SU(2)$.
\efolg
This corollary together with corollary \ref{folg1} proves proposition
\ref{satz1} and therefore theorem \ref{theo1}.

\section{Consequences}
In order to explain the conclusion more in detail we cite the result of L. Berard-Bergery  and A. Ikemakhen
about four different types of indecomposable, non-irreducible Lorentzian holonomy algebras. One considers $\rrn$ with
the Minkowskian scalar product of the form
$\eta := \left( \begin{array}{ccc}
0	 &   0^t&1	\\
0	 & E_n & 0 \\
1	 & 0^t &0
\end{array}\right)$. Any indecomposably, but non-irreducibly acting subalgebra of $\laso(\eta)$ is contained in the
parabolic algebra $(\rr\+\lason)\ltimes \rrn$. Furthermore one can prove the following result.
\btheo\label{theoI}
 \cite{bb-ike93}
Let ${\frak h}$ be a subalgebra of $\laso(\eta)$ which acts indecomposably and non-irreducibly
on $\rr^{n+2}$, $\lag:= pr_{\lason}(\mf{h})$ with $\lag=\mf{z}\+\mf{d}$ its Levi-decomposition in the center and the derived Lie algebra.
Then $\mf{h}$ belongs to one of the following types.
\begin{enumerate}
\item If ${\frak h}$ contains $\rrn$, then we have the types
\begin{description}
\item[ Type 1:] $ {\frak h}$ contains $\rr$.
Then $ {\frak h} =\left(\rr \oplus {\frak g} \right) \ltimes \rrn  $.
\item[ Type 2:]  $ pr_{\rr}({\frak h})= 0$ i.e. ${\frak h} ={\frak g}\ltimes
\rrn$.
\item[ Type 3:] Neither Type 1 nor Type 2.

In that case there exists a surjective homomorphism $\varphi: {\frak z} \rightarrow
\rr $, such that
\[{\frak h} = \left( {\frak l} \oplus \mf{d} \right)\ltimes \rrn\]
where ${\frak l}:= graph\ \varphi = \{ ( \varphi(T),T) | T\in {\frak z} \} \subset
\rr\oplus{\frak z}  $.
Or written as matrices:
\[{\frak h}= \left\{ \left.
\left(
\begin{array}{ccc}
 \varphi(A) 	& 	v^t	& 0\\
  0	&A+B& -v\\
0&0& -\varphi(A)\\
\end{array}
\right) \right| A\in {\frak z} , B\in \mf{d}, v\in \mathbb{R}^n\right\}.
\]
\end{description}
\item In case ${\frak h}$ does
not contain $\rrn$ we have {\em {\bf Type 4:}}

There exists
\begin{enumerate}
\item a {non-trivial} decomposition $\rrn =
\mathbb{R}^k\oplus \mathbb{R}^l$, $0<k,l<n$,
\item a surjective homomorphism $\varphi: {\frak z} \rightarrow \rr^l$
\end{enumerate}
such that ${\frak g}\subset {\frak s}{\frak o}(k)$ and
$ {\frak h} =  \left( \mf{d} \oplus {\frak l} \right)\ltimes \rr^k \subset
{\frak p}$
where ${\frak l}:= \{ \left( \varphi(T), T \right) | T\in {\frak z} \} = graph\
\varphi \subset \rr^l \oplus{\frak z}$.
Or written as matrices:
\[{\frak h}= \left\{ \left.
\left(
\begin{array}{cccc}
0	 & \varphi(A)^t 	& 	v^t	& 0\\
0	 &  0	&A+B& -v\\
0&0&0& -\varphi(A)\\
0&0&0&0
\end{array}
\right) \right| A\in {\frak z} , B\in \mf{d}, v\in \mathbb{R}^k \right\}.
\]
\end{enumerate}
\end{theo}
From this and theorem \ref{theo1} we get an obvious corollary.
\bfolg
Let $\mf{h}$ be the holonomy algebra of an indecomposable, non irreducible $n+2$-dimensional Lorentzian manifold.
\bnum
\item
If $\mf{h}$ is of type 1 or 2, then it holds $\mf{h}=(\rr\+\lag)\ltimes \rrn$ or
$\lag\ltimes \rrn$, where $\lag$ is a Riemannian holonomy algebra.
\item
If $\mf{h}$ is of type 3 or 4, then $\lag=pr_{\lason}\mf{h}$ is a Riemannian holonomy algebra
with at least one irreducible factor equal to a Riemannian holonomy algebra with center, i.e. equal to
$\laso(2)$ acting on $\rr^2$ or on itself, $\laso(2)\+\laso(n)$ acting on $\rr^{2n}$
$\laso(2)\+\laso(10)$ acting on $\rr^{32}$ as
the reellification of the complex spinor module of dimension 16,  $\laso(2)\+ \mf{e}_6$ acting on $\rr^{54}$,
$\mf{u}(n)$ acting on $\rr^{2n}$ or on $\rr^{n(n-1)}$.
\enum
\efolg
Regarding the epimorphisms $\vf:\mf{z}(\lag)\mapsto \rr^k$ from theorem \ref{theoI} there is a theorem of C. Boubel
\cite[Th\'{e}or\`{e}me 3.IV.3 and Corollaire 3.IV.3]{boubel00} which describes how to construct a
metric with holonomy of type 3 or 4 from metrics with holonomy of type 1 or 2, under certain algebraic conditions
on $\lag$ of course. Knowing the possible $\lag$'s gives some candidates to start with in order
to construct such metrics of type 3 and 4.

\bigskip

Finally we want to draw some conclusions about the existence of parallel spinor fields on Lorentzian manifold.
The existence of a parallel spinor field on a Lorentzian spin manifold implies the existence of
a parallel vector field which has to be lightlike or timelike. In the latter case the manifold splits by the
de-Rham decomposition theorem (at least locally) into
a factor $(\rr,-dt^2)$ and Riemannian factors which are flat or irreducible with a parallel spinor, i.e. with holonomy
$\{1\},\ G_2,\ Spin(7),\ Sp(k) $ or $SU(k)$.

In the case where the parallel vector field is lightlike  we have a Lorentzian factor which is indecomposable, but with
parallel lightlike vector field (and parallel spinor) and flat or irreducible Riemannian manifolds with parallel spinors.
Hence in this case one has to know which indecomposable Lorentzian manifolds admit a parallel spinor. The existence of the
lightlike parallel vector field forces the holonomy
of such a manifold with parallel spinor to be contained in $\lason\ltimes \rrn$ i.e. to be of type 2 or 4. Furthermore the
 spin representation of $\lason$-projection $\lag\subset\lason$ must admit a trivial subrepresentation (see for example
 \cite{leistner01}). Up to dimension $n+2=11$ these groups where described by R. L. Bryant in \cite{bryant00} and
 J.M. Figueroa O'Farrill in \cite{farrill99mw}
obtaining that at least the maximal ones are of type 2 and of the shape $(Riemannian\ holonomy)\ltimes\rrn$. Now we get this
result in general.
\bfolg
Let $\mf{h}$ be the holonomy algebra of an indecomposable Lorentzian spin manifold with parallel spinor field.
Then $\mf{h}=\lag\ltimes \rrn$ where $\lag$ is the holonomy algebra of a Riemannian manifold with parallel spinor,
i.e. a sum of the following algebras: $\{0\},\ \mf{g}_2,\ \mf{spin}(7),\ \mf{sp}(k) $ or $\mf{su}(k)$.
\efolg
\bprf
We have to exclude that the holonomy algebra can be of type 4 under the assumption of a parallel spinor field.
But if $\mf{h}$ is of type 4, the $\lason$-projection $\lag$ has a $\laso(2)$-summand. (Also in the $\mf{u}(n)$ case since
$\mf{u}(n)=\laso(2)\+\mf{su}(n)$.) But a direct calculation ($\laso(2)=\rr J$ with $J^2=-id$)
shows that the spin representation of such a $\laso(2)$-summand
is an isomorphism of the spinor module, i.e. there can be no trivial subrepresentation.
\eprf
Finally one should remark that it is very desirable to find a direct proof of
these facts avoiding this cumbersome case-by-case analysis.

\bibliography{GEOBIB,SPINBIB,HOLBIB,ALGBIB,thomas}
{\small
{\tt leistner@mathematik.hu-berlin.de}\\
{\sc Institut f\"{u}r Mathematik, Humboldt-Universit\" {a}t Berlin\\
Unter den Linden 6, D-10099 Berlin}}

\end{document}